\documentclass{elsart}

\usepackage{amssymb}
\usepackage[latin1]{inputenc}

\begin{document}

\begin{frontmatter}



\title{Accurate computations with Said-Ball-Vandermonde matrices}


\author{Ana Marco
\corauthref{1}}
\corauth[1]{Corresponding author.}
\ead{ana.marco@uah.es},
\author{José-Javier Mart{\'\i}nez}
\ead{jjavier.martinez@uah.es}

%\author{A. Marco},
%\ead{ana.marco@uah.es}
%\author{J. J. Mart{\'\i}nez
%\corauthref{1}}
%\ead{jjavier.martinez@uah.es}
%\corauth[1]{Corresponding author.}

\address{Departamento de Matemáticas, Universidad de Alcalá, Campus Universitario, 28871 Alcalá de Henares, Madrid, Spain}

\begin{abstract}

A generalization of the Vandermonde matrices which arise when the power basis is replaced by the Said-Ball basis is considered. When the nodes are inside the interval $(0,1)$, then those matrices are strictly totally positive. An algorithm for computing the bidiagonal decomposition of those Said-Ball-Vandermonde matrices is presented, which allows to use known algorithms for totally positive matrices represented by their bidiagonal decomposition. The algorithm is shown to be fast and to guarantee high relative accuracy. Some numerical experiments which illustrate the good behaviour of the algorithm are included.


\end{abstract}


\end{frontmatter}

\section{Introduction}

Numerical computing with structured totally nonnegative matrices
is a classical subject in the field of numerical linear algebra
which has recently received a renewed attention, as can be seen in
the recent survey paper \cite{DDHK}, where several different classes of
structured matrices are considered, among them totally positive matrices.

Classically, a matrix is said to be {\it totally positive} if all its minors are nonnegative \cite{GP96}. Consequently, the
 matrices with that property are also called {\it totally nonnegative matrices} \cite{FALL}, and this term is becoming more used in recent literature.

The fact that a nonsingular totally nonnegative (TN) matrix can be
decomposed as a product of nonnegative bidiagonal factors was used
by Koev \cite{KOEV05, KOEV07} to develop several accurate algorithms for the
general class of TN matrices. A detailed survey of several
results related to TN matrices, including the bidiagonal
factorization, has been presented in \cite{FALL}.

Nevertheless, it must be stressed that the algorithms of Koev
\cite{KOEV05, KOEV07} start from the bidiagonal decomposition of a TN matrix
$A$, which is stored in a matrix which is denoted there as $\mathcal{BD}(A)$, and that such a
decomposition needs to be computed for each {\it particular} class
of TN matrices being considered. Using  the words of the section
devoted to conclusions and open problems in \cite{KOEV07}:

{\it The caveat in our algorithms is that every TN matrix must be
represented by its bidiagonal decomposition. While every TN matrix
intrinsically possesses such a decomposition, and for many classes
of structured matrices this decomposition is very easy to obtain
accurately , there are important TN matrices for which we know of
no accurate and efficient way to compute their bidiagonal
decompositions.}

Examples of totally nonnegative matrices for which there are accurate and efficient algorithms for computing $\mathcal{BD}(A)$ are Vandermonde
\cite{BP, GV, HIG87, HIG02}, Cauchy \cite{BKO}, Cauchy-Vandermonde \cite{MP98, MP03}, generalized
Vandermonde \cite{DK} and Bernstein-Vandermonde matrices \cite{MM07}.

On the other hand, it is not always recognized that while {\it
Neville elimination} \cite{FALL, GP92, GP94, GP96} is a key theoretical tool for
the analysis of that bidiagonal decomposition, it generally fails
to provide an accurate algorithm for computing $\mathcal{BD}(A)$. This
fact is explicitly noted in \cite{KOEV}, where the author
indicates that {\it the function TNBD is the only function in the
package TNTool that does not guarantee high relative accuracy}.

Consequently, the accurate (and, if possible, fast) computation of
$\mathcal{BD}(A)$ is a previous task to be performed before applying Koev's
algorithms to a given class of TN matrices. The importance of those algorithms was
very acknowledged in \cite{WAT}, while relevant previous results were presented in \cite{DetAl}.

In this work we are extending to the class of Said-Ball-Vandermonde matrices the work we
have recently carried out for the class of Bernstein-Vandermonde
matrices \cite{MM07}. A crucial fact for obtaining high relative accuracy in our algorithm is that it
satisfies what is called in \cite{DDHK} the NIC (no inaccurate cancellation) condition:

{\bf NIC:} The algorithm only multiplies, divides, adds (resp., substracts) real numbers with like (resp., differing) signs, and otherwise only adds or substracts input data.

The Said-Ball basis is a generalization of the Ball basis \cite{BALL1, BALL2, BALL3}, a
well-known basis for cubic polynomials on a finite interval which is useful in the field of Computer-Aided
Design. The Said-Ball basis was introduced for odd degree polynomials by Said in \cite{SAID}, and then its definition for  polynomials of even degree was suggested in \cite{HU}. Its properties in connection with
total positivity and shape preservation were studied by Goodman
and Said  for odd degree polynomials \cite{GS}, and recently by Delgado an Pe\~na in \cite{DP06}, where it was
established that the Said-Ball basis is a normalized totally
positive (NTP) basis for every value of the polynomial degree.

The rest of the paper is organized as follows. Some basic results on Neville elimination
and total positivity are recalled in Section 2. In Section 3 the
bidiagonal decomposition of a Said-Ball-Vandermonde
matrix and of its inverse are presented. The algorithm for computing these bidiagonal factorizations is introduced in Section 4. In Section 5 the problems of linear system solving and eigenvalue computation for a Said-Ball-Vandermonde
matrix are considered. Finally, Section 6 is devoted to illustrate the accuracy of our algorithms by means of some
numerical experiments.

\section{Basic results on Neville elimination and total positivity} %%%%%%%%%%%%%%%%%%%%%%%%%%%%%%%%%%%%%%%%%%%%%%%%%%%%%%%%%%%%%%%%

In this section we will briefly recall some basic results on Neville
elimination and total positivity which we will apply in Section 3.
Our notation follows the notation used in \cite{GP92, GP94}. Given $k$,
$n \in {\bf N}$  ($1 \leq k \leq n$), $Q_{k,n}$ will denote the set
of all increasing sequences of $k$ positive integers less than or
equal to $n$.

Let $A$ be a real square matrix of order $n$. For $k \leq n$, $m
\leq n$, and for any $\alpha \in Q_{k,n}$ and $\beta \in Q_{m,n}$,
we will denote by $A[\alpha \vert \beta ]$ the submatrix $k\times m$
of $A$ containing the rows numbered by $\alpha $ and the columns
numbered by $\beta $.

The fundamental tool for obtaining the theoretical results applied in this
paper is the {\it Neville elimination} (see \cite{GP92, GP94}), a procedure that
makes zeros in a matrix adding to a given row an appropriate
multiple of the previous one. For a nonsingular matrix
$A=(a_{i,j})_{1\leq i,j\leq n}$, it consists on $n-1$ steps
resulting in a sequence of matrices $A:=A_1\to A_2\to \ldots \to
A_n$, where $A_t= (a_{i,j}^{(t)})_{1\leq i,j\leq n}$ has zeros below
its main diagonal in the $t-1$ first columns. The matrix $A_{t+1}$
is obtained from $A_t$ ($t=1,\ldots ,n$) by using the following
formula:
$$
a_{i,j}^{(t+1)}:= \left\{ \begin{array}{ll} a_{i,j}^{(t)}~, &
\textnormal{if} \quad i\leq t\\
a_{i,j}^{(t)}-(a_{i,t}^{(t)}/a_{i-1,t}^{t})a_{i-1,j}^{(t)}~,~ &
\textnormal{if} \quad i\geq t+1 ~\textnormal{and}~ j\geq t+1\\
0~, & \textnormal{otherwise}.
\end{array}\right.  \eqno(2.1)
$$
In this process the element
$$
p_{i,j}:=a_{i,j}^{(j)} \qquad 1\leq j\leq n; ~~ j\leq i\leq n
$$
is called {\it pivot} ($i,j$) of the Neville elimination of $A$. The
process would break down if any of the pivots $p_{i,j}$ ($j\leq
i<n$) is zero. In that case we can move the corresponding rows to
the bottom and proceed with the new matrix, as described in \cite{GP92}.
The Neville elimination can be done without row exchanges if all the
pivots are nonzero, as it will happen in our situation. The pivots
$p_{i,i}$ are called {\it diagonal pivots}. If all the pivots
$p_{i,j}$ are nonzero, then $p_{i,1}=a_{i,1}\, \forall i$ and, by
Lemma 2.6 of \cite{GP92}
$$p_{i,j}={\det A[i-j+1,\ldots ,i\vert 1,\ldots ,j]\over \det
A[i-j+1,\ldots ,i-1\vert 1,\ldots ,j-1]}\qquad 1<j\leq i\leq n.
\eqno(2.2)$$

The element
$$
m_{i,j}=\frac{p_{i,j}}{p_{i-1,j}} \qquad 1\leq j\leq n; ~~ j< i\leq
n  \eqno(2.3)
$$
is called {\it multiplier} of the Neville elimination of $A$. The
matrix $U:=A_n$ is upper triangular and has the diagonal pivots in
its main diagonal.

The {\it complete Neville elimination} of a matrix $A$ consists on
performing the Neville elimination of $A$ for obtaining $U$ and then
continue with the Neville elimination of $U^T$. The pivot
(respectively, multiplier) $(i,j)$ of the complete Neville
elimination of $A$ is the pivot (respectively, multiplier) $(j,i)$
of the Neville elimination of $U^T$, if $j\ge i$. When no row
exchanges are needed in the Neville elimination of $A$ and $U^T$, we
say that the complete Neville elimination of $A$ can be done without
row and column exchanges, and in this case the multipliers of the
complete Neville elimination of $A$ are the multipliers of the
Neville elimination of $A$ if $i\ge j$ and the multipliers of the
Neville elimination of $A^T$ if $j\ge i$.

A matrix is called {\it strictly totally positive}) if all its minors are positive. The Neville elimination characterizes the
strictly totally positive matrices as follows (\cite{GP92}):

{\bf Theorem 2.1.} A matrix is strictly totally positive if and only
if its complete Neville elimination can be performed without row and
column exchanges, the multipliers of the Neville elimination of $A$
and $A^T$ are positive, and the diagonal pivots of the Neville
elimination of $A$ are positive.

As it can be seen in \cite{DP06}, the Said-Ball-Vandermonde matrices are strictly totally positive when the real numbers
satisfy $0<t_1<t_2<\ldots<t_{n+1}<1$, and this fact has inspired our
search for a fast algorithm, but this result will also be shown to
be a consequence of our Theorem 3.2.

\section{Bidiagonal decomposition} %%%%%%%%%%%%%%%%%%%%%%%%%%%%%%%%%%%%%%%%%%%%%%%%%%%%%%%%%%%%%%%%%%%%%%%%%%%%%%%%%%%%%%%%%%%%%%%%%

The {\it Said-Ball basis} $\mathcal{S}_n=\{s_0^n(t), s_1^n(t), \ldots, s_n^n(t) \}$ of the space $\Pi_n(t)$ of the polynomials of degree less than or equal to $n$ on the interval $[0, 1]$ is defined by:
$$
\begin{array}{ll}
s_i^n(t)={\lfloor n/2 \rfloor +i \choose i}t^i(1-t)^{\lfloor n/2 \rfloor +1}, & \quad 0\leq i \leq \lfloor (n-1)/2 \rfloor,\\
s_i^n(t)={ \lfloor n/2 \rfloor +n-i\choose n-i}t^{\lfloor n/2 \rfloor +1}(1-t)^{n-i}, & \quad \lfloor n/2 \rfloor+1 \leq i \leq n,
\end{array}
$$
and, if $n$ is even
$$
s_{n/2}^n(t)={n \choose n/2}t^{n/2}(1-t)^{n/2},
$$
where $\lfloor m \rfloor$ is the greatest integer less than or equal to $m$.

From now on, we will call {\it Said--Ball--Vandermonde matrices} (SB--Vandermonde matrices in the sequel) the generalization of the Vandermonde matrices obtained when considering the Said-Ball basis instead of the power basis. The SB--Vandermonde matrices are therefore
$$A=
\left(
    \begin{array}{cccc}
      {\frac{n-1}{2} \choose 0}(1-t_1)^{\frac{n+1}{2}} & {\frac{n-1}{2} \choose 0}(1-t_2)^{\frac{n+1}{2}} & \cdots & {\frac{n-1}{2} \choose 0}(1-t_{n+1})^{\frac{n+1}{2}}\\
      {\frac{n-1}{2}+1 \choose 1}t_1(1-t_1)^{\frac{n+1}{2}} & {\frac{n-1}{2}+1 \choose 1}t_2(1-t_2)^{\frac{n+1}{2}} & \cdots & {\frac{n-1}{2}+1 \choose 1}t_{n+1}(1-t_{n+1})^{\frac{n+1}{2}}\\
      \vdots & \vdots & \ddots & \vdots \\
      {n-1 \choose \frac{n-1}{2}}t_1^{\frac{n-1}{2}}(1-t_1)^{\frac{n+1}{2}} & {n-1 \choose \frac{n-1}{2}}t_2^{\frac{n-1}{2}}(1-t_2)^{\frac{n+1}{2}} & \cdots & {n-1 \choose \frac{n-1}{2}}t_{n+1}^{\frac{n-1}{2}}(1-t_{n+1})^{\frac{n+1}{2}}\\
      {n-1 \choose \frac{n-1}{2}}t_1^{\frac{n+1}{2}}(1-t_1)^{\frac{n-1}{2}} & {n-1 \choose \frac{n-1}{2}}t_2^{\frac{n+1}{2}}(1-t_2)^{\frac{n-1}{2}} & \cdots & {n-1 \choose \frac{n-1}{2}}t_{n+1}^{\frac{n+1}{2}}(1-t_{n+1})^{\frac{n-1}{2}}\\
      \vdots & \vdots & \ddots & \vdots\\
      {\frac{n-1}{2}+1 \choose 1}t_1^{\frac{n+1}{2}}(1-t_1) & {\frac{n-1}{2}+1 \choose 1}t_2^{\frac{n+1}{2}}(1-t_2) & \cdots & {\frac{n-1}{2}+1 \choose 1}t_{n+1}^{\frac{n+1}{2}}(1-t_{n+1})\\
      {\frac{n-1}{2} \choose 0}t_1^{\frac{n+1}{2}} & {\frac{n-1}{2} \choose 0}t_2^{\frac{n+1}{2}} & \cdots & {\frac{n-1}{2} \choose 0}t_{n+1}^{\frac{n+1}{2}}
    \end{array}
  \right)^T
$$
in the case of odd $n$, and
$$A=
\left(
    \begin{array}{cccc}
      {\frac{n}{2} \choose 0}(1-t_1)^{\frac{n+2}{2}} & {\frac{n}{2} \choose 0}(1-t_2)^{\frac{n+2}{2}} & \cdots & {\frac{n}{2} \choose 0}(1-t_{n+1})^{\frac{n+2}{2}}\\
      {\frac{n}{2}+1 \choose 1}t_1(1-t_1)^{\frac{n+2}{2}} & {\frac{n}{2}+1 \choose 1}t_2(1-t_2)^{\frac{n+2}{2}} & \cdots & {\frac{n}{2}+1 \choose 1}t_{n+1}(1-t_{n+1})^{\frac{n+2}{2}}\\
      \vdots & \vdots & \ddots & \vdots \\
      {n-1 \choose \frac{n-2}{2}}t_1^{\frac{n-2}{2}}(1-t_1)^{\frac{n+2}{2}} & {n-1 \choose \frac{n-2}{2}}t_2^{\frac{n-2}{2}}(1-t_2)^{\frac{n+2}{2}} & \cdots &  {n-1 \choose \frac{n-2}{2}}t_{n+1}^{\frac{n-2}{2}}(1-t_{n+1})^{\frac{n+2}{2}}\\
      {n \choose \frac{n}{2}}t_1^{\frac{n}{2}}(1-t_1)^{\frac{n}{2}} & {n \choose \frac{n}{2}}t_2^{\frac{n}{2}}(1-t_2)^{\frac{n}{2}} & \cdots & {n \choose \frac{n}{2}}t_{n+1}^{\frac{n}{2}}(1-t_{n+1})^{\frac{n}{2}}\\
      {n-1 \choose \frac{n-2}{2}}t_1^{\frac{n+2}{2}}(1-t_1)^{\frac{n-2}{2}} & {n-1 \choose \frac{n-2}{2}}t_2^{\frac{n+2}{2}}(1-t_2)^{\frac{n-2}{2}} & \cdots & {n-1 \choose \frac{n-2}{2}}t_{n+1}^{\frac{n+2}{2}}(1-t_{n+1})^{\frac{n-2}{2}}\\
      \vdots & \vdots & \ddots & \vdots \\
      {\frac{n}{2}+1 \choose 1}t_1^{\frac{n+2}{2}}(1-t_1) & {\frac{n}{2}+1 \choose 1}t_2^{\frac{n+2}{2}}(1-t_2) & \cdots & {\frac{n}{2}+1 \choose 1}t_{n+1}^{\frac{n+2}{2}}(1-t_{n+1})\\
      {\frac{n}{2} \choose 0}t_1^{\frac{n+2}{2}} &  {\frac{n}{2} \choose 0}t_2^{\frac{n+2}{2}} & \cdots &  {\frac{n}{2} \choose 0}t_{n+1}^{\frac{n+2}{2}}

\end{array}
  \right)^T
$$
in the case of even $n$.

 It must be observed that the SB--Vandermonde matrix $A$ is the coefficient matrix associated with the following interpolation problem in the Said-Ball basis $\mathcal{S}_n$: given the interpolation nodes $\{t_i: ~i=1, \ldots, n+1\}$ and the interpolation data $\{b_i: ~i=1, \ldots, n+1\}$ find the polynomial
$$
p(t)=\sum_{k=0}^n a_k s_k^n(t)
$$
such that $p(t_i)=b_i$ for $i=1,\ldots,n+1$.

From now on, we will assume $0<t_1<t_2< \ldots<t_{n+1}<1$.

\medskip
{\bf Proposition 3.1.} {\it The determinant of the SB--Vandermonde matrix $A$ defined above is
$$
\det A=\Bigg[ {\frac{n-1}{2} \choose 0} {\frac{n+1}{2} \choose 1} {\frac{n+3}{2} \choose 2} \cdots {n-2 \choose \frac{n-3}{2}} {n-1 \choose \frac{n-1}{2}} \Bigg]^2 \prod_{1 \leq i < j \leq n+1}(t_j-t_i),
$$
if $n$ is odd, and
$$
\det A=\Bigg[ {\frac{n}{2} \choose 0} {\frac{n+2}{2} \choose 1} {\frac{n+4}{2} \choose 2} \cdots {n-1 \choose \frac{n-2}{2}}\Bigg]^2 {n \choose \frac{n}{2}} \prod_{1 \leq i < j \leq n+1}(t_j-t_i),
$$
if $n$ is even.}

{\bf Proof.} Here we include the proof for the case in which $n$ is odd. The proof in the even case is completely analogous.

Looking at \cite{HU}, it can be observed that the matrix of change of basis from the Bernstein basis to the Said-Ball basis is a block-diagonal matrix $M$ with triangular diagonal blocks, and whose determinant is
$$
\det M=\frac{\Big[{\frac{n-1}{2} \choose 0}{\frac{n+1}{2} \choose 1}{\frac{n+3}{2} \choose 2}\cdots {n-2 \choose \frac{n-3}{2}}{n-1 \choose \frac{n-1}{2}}\Big]^2}{{n \choose 0}{n \choose 1} \cdots {n \choose n}}. \eqno(3.1)
$$

As it can be seen, for example, in \cite{MM07}, the matrix of change of basis from the Bernstein basis
$$
\mathcal{B}_n=\big\{ b_i^{(n)}(t) = {n \choose i} (1 - t)^{n-i} t^i,
\qquad i = 0, \ldots, n \big\}
$$
to the power basis $\{1, t, t^2, \ldots, t^n \}$ is a lower triangular matrix $N$ of order $n+1$ whose determinant is
$$
\det N= {n \choose 0}{n \choose 1} \cdots {n \choose n}. \eqno(3.2)
$$
Taking this into account, the matrix of change of basis from the power basis to the Said-Ball basis is $MN^{-1}$, and consequently,
$$
\det A=\frac{\det M}{\det N}\det V ,
$$
where $V$ is the Vandermonde matrix
$$
V=\left( \begin{array}{ccccc}
1 & t_1 & t_1^2 & \cdots & t_1^n\\
1 & t_2 & t_2^2 & \cdots & t_2^n\\
\vdots & \vdots & \ddots & \vdots\\
1 & t_{n+1} & t_{n+1}^2 & \cdots & t_{n+1}^n
\end{array} \right).
$$
Using the well-known formula for the determinant of  a Vandermonde
matrix
$$
\det V=\prod_{1\leq i < j \leq n+1}(t_j-t_i)
$$
and the equations (3.1) and (3.2), the proof is concluded.
$\Box$

The following two theorems will be essential in the construction of our algorithm for computing $\mathcal{BD}(A)$ of a SB--Vandermonde matrix.

\medskip
{\bf Theorem 3.2.} Let $A=(a_{i,j})_{1\le i,j\le n+1}$ be a
SB--Vandermonde matrix whose nodes satisfy $0 < t_1 < t_2
<\ldots < t_n <t_{n+1} <1$. Then $A^{-1}$ admits a factorization in
the form
$$A^{-1}=G_1G_2\cdots G_{n}D^{-1}F_{n}F_{n-1}\cdots F_1, \eqno(3.3)$$
where $G_i$ are upper triangular bidiagonal matrices, $F_i$ are
lower triangular bidiagonal matrices ($i=1,\ldots,n$), and $D$ is a
diagonal matrix.

{\bf Proof.} The matrix $A$ is a strictly totally positive matrix (see \cite{CP,DP06}) and therefore, by Theorem 2.1, the complete Neville elimination of
$A$ can be performed without row and column exchanges providing the
following factorization of $A^{-1}$ (see \cite{GP92, GP94}):
$$A^{-1}=G_1G_2\cdots G_{n}D^{-1}F_{n}F_{n-1}\cdots F_1, $$
where $F_i$ ($1\le i\le n$) are bidiagonal matrices of the form
$$F_i=\left( \begin{array}{cccccccc}
1 & & & & & & & \\
0 & 1 & & & & & & \\
& \ddots & \ddots & & & & & \\
& & 0 & 1 & & & & \\
& & & -m_{i+1,i} & 1 & & & \\
& & & & -m_{i+2,i} & 1 & & \\
& & & & & \ddots & \ddots & \\
& & & & & & -m_{n+1,i} & 1
\end{array} \right),  \eqno(3.4)
$$
$G^T_i$ ($1\le i\le n$) are bidiagonal matrices of the form
$$G_i^T=\left( \begin{array}{cccccccc}
1 & & & & & & & \\
0 & 1 & & & & & & \\
& \ddots & \ddots & & & & & \\
& & 0 & 1 & & & & \\
& & & -\widetilde m_{i+1,i} & 1 & & & \\
& & & & -\widetilde m_{i+2,i} & 1 & & \\
& & & & & \ddots & \ddots & \\
& & & & & & -\widetilde m_{n+1,i} & 1
\end{array} \right),  \eqno(3.5)
$$
and $D$ is the diagonal matrix whose $i$th ($1\le i\le n+1$)
diagonal entry is the diagonal pivot $p_{i,i}=a_{i,i}^{(i)}$ of the
Neville elimination of $A$:
$$
D=\textnormal{diag}\{p_{1,1},p_{2,2},\ldots,p_{n+1,n+1}\}. \eqno(3.6)
$$

First we obtain the expressions for the multipliers $m_{i,j}$ and $\widetilde m_{i,j}$, and for the diagonal pivots $p_{i,i}$ in the case of odd $n$.

Taking into account that the minors of $A$ with $j$ initial
consecutive columns and $j$ consecutive rows starting with row $i$
are
$$
\begin{array}{l}
\det A[i, \ldots, i+j-1 \vert 1, \ldots, j] = {\frac{n-1}{2} \choose 0} {\frac{n-1}{2}+1 \choose 1}\cdots {\frac{n-1}{2}+j-1 \choose j-1}\\  (1-t_{i})^{\frac{n+1}{2}} (1-t_{i+1})^{\frac{n+1}{2}}\cdots(1-t_{i+j-1})^{\frac{n+1}{2}} \prod_{i \leq k < l \leq i+j-1}(t_l-t_k),
\end{array}
$$
if $j \leq \frac{n+1}{2}$, and
$$
\begin{array}{l}
\det A[i, \ldots, i+j-1 \vert 1, \ldots, j] = {\frac{n-1}{2} \choose 0} {\frac{n-1}{2}+1 \choose 1}\cdots {\frac{n-1}{2}+n-j \choose n-j}\\
\Big[ {\frac{n-1}{2}+n-j+1  \choose n-j+1} \cdots {n-1 \choose \frac{n-1}{2}} \Big]^2
(1-t_{i})^{n-j+1} (1-t_{i+1})^{n-j+1}\cdots(1-t_{i+j-1})^{n-j+1}\\
 \prod_{i \leq k < l \leq i+j-1}(t_l-t_k)
\end{array}
$$
if $j>\frac{n+1}{2}$,

a result that follows from the properties of the determinants and Proposition 3.1, and that $m_{i,j}$ are the multipliers of the Neville elimination of $A$, we obtain that
$$
m_{i,j}=\left\{ \begin{array}{ll}
\frac{(1-t_i)^{\frac{n+1}{2}} \prod_{k=1}^{j-1}(t_i-t_{i-k})}{(1-t_{i-1})^{\frac{n+1}{2}}\prod_{k=2}^{j}(t_{i-1}-t_{i-k})}, & j=1,\ldots,\frac{n+1}{2}; ~i=j+1, \ldots,n+1,\\
 & \\
\frac{(1-t_i)^{n-j+1} (1-t_{i-j})\prod_{k=1}^{j-1}(t_i-t_{i-k})}{(1-t_{i-1})^{n-j+2}\prod_{k=2}^{j}(t_{i-1}-t_{i-k})}, & j=\frac{n+3}{2},\ldots,n; ~i=j+1,\ldots, n+1.
\end{array} \right. \eqno(3.7)
$$
As for the minors of $A^T$ with $j$ initial consecutive columns and $j$ consecutive rows starting with row $i$, they are:
$$
\begin{array}{l}
\det A^T[i, \ldots, i+j-1 \vert 1,\ldots,j]=
{\frac{n-1}{2} +i-1 \choose i-1}{\frac{n-1}{2} +i \choose i}\cdots {\frac{n-1}{2}+i+j-2 \choose i+j-2}\\
(1-t_1)^{\frac{n+1}{2}}(1-t_2)^{\frac{n+1}{2}}\cdots (1-t_j)^{\frac{n+1}{2}}t_1^{i-1}t_2^{i-1}\cdots t_j^{i-1} \prod_{1 \leq k <l \leq j}(t_l-t_k),
\end{array}
$$
if $i\leq \frac{n+1}{2}$ and $i+j-1\leq \frac{n+1}{2}$,
$$
\begin{array}{l}
\det A^T[i, \ldots, i+j-1 \vert 1,\ldots,j]={\frac{n-1}{2} +i-1 \choose i-1}{\frac{n-1}{2} +i \choose i}\cdots {n-1 \choose \frac{n-1}{2}}
{n-1 \choose \frac{n-1}{2}}\\
{ n-2 \choose \frac{n-3}{2}}\cdots {\frac{n-1}{2}+n-i-j+2 \choose n-i-j+2}
t_1^{i-1}t_2^{i-1}\cdots t_j^{i-1} (1-t_1)^{n-i-j+2}\\
(1-t_2)^{n-i-j+2}\cdots (1-t_j)^{n-i-j+2}\prod_{1 \leq k <l \leq j}(t_l-t_k),
\end{array}
$$
if $i\leq \frac{n+1}{2}$ and $i+j-1> \frac{n+1}{2}$, and
$$
\begin{array}{l}
\det A^T[i, \ldots, i+j-1 \vert 1,\ldots,j]=
{\frac{n-1}{2}+n-i+1 \choose n-i+1}{\frac{n-1}{2}+n-i \choose n-i}\cdots {\frac{n-1}{2}+n-i-j+2 \choose n-i-j+2}\\
t_1^{\frac{n+1}{2}}t_2^{\frac{n+1}{2}} \cdots t_j^{\frac{n+1}{2}}(1-t_1)^{n-i-j+2}
(1-t_2)^{n-i-j+2}\cdots (1-t_j)^{n-i-j+2}\\
\prod_{1 \leq k <l \leq j}(t_l-t_k),
\end{array}
$$
if $i> \frac{n+1}{2}$.

These expressions also follow from the properties of the determinants
and Proposition 3.1. Since the entries $\widetilde m_{i,j}$ are the
multipliers of the Neville elimination of $A^T$, using the previous
expressions for the minors of $A^T$ with initial consecutive columns
and consecutive rows, it is obtained that
$$
\widetilde m_{i,j}=\left\{ \begin{array}{ll}
\frac{\frac{n-1}{2}+i-1}{i-1}t_j, & \quad i=2,\ldots,\frac{n+1}{2}; ~j=1,\ldots,i-1,\\
\frac{t_j}{\prod_{k=1}^j (1-t_k)} & \quad i=\frac{n+3}{2}; ~j=1,\ldots,\frac{n+1}{2},\\
\frac{n-i+2}{\frac{n-1}{2}+n-i+2}\frac{1}{1-t_j} & \quad i=\frac{n+5}{2},\ldots,n+1; j=1,\ldots,i-\frac{n+3}{2},\\
\frac{n-i+2}{\frac{n-1}{2}+n-i+2}\frac{t_j}{1-t_j} & \quad i=\frac{n+5}{2},\ldots,n+1; j=i-\frac{n+1}{2},\ldots,i-1.
\end{array} \right. \eqno(3.8)
$$
Finally, the diagonal entries of $D$ are:
$$
p_{i,i}=\left\{ \begin{array}{ll}
{\frac{n-1}{2}+i-1 \choose i-1} (1-t_i)^{\frac{n+1}{2}}\prod_{k<i}(t_i-t_k), & \quad i=1,\ldots,\frac{n+1}{2},\\
 & \\
{\frac{n-1}{2}+n-i+1\choose n-i+1} \frac{(1-t_i)^{n-i+1} \prod_{k<i}(t_i-t_k)}{\prod_{k=1}^{i-1}(1-t_k)}, & \quad i=\frac{n+3}{2}, \ldots, n+1.
\end{array} \right. \eqno(3.9)
$$
The formulas for $p_{i,i}$ are obtained by using the expressions for the minors of $A$ with initial consecutive columns and initial consecutive rows.

As for the case in which $n$ is even, proceeding analogously as in the odd case we obtain the following expressions for the multipliers $m_{i,j}$ and $\widetilde m_{i,j}$, and the diagonal pivots $p_{i,i}$:
$$
m_{i,j}=\left\{ \begin{array}{ll}
\frac{(1-t_i)^{\frac{n+2}{2}} \prod_{k=1}^{j-1}(t_i-t_{i-k})}{(1-t_{i-1})^{\frac{n+2}{2}}\prod_{k=2}^{j}(t_{i-1}-t_{i-k})}, & j=1,\ldots,\frac{n}{2}; ~i=j+1, \ldots,n+1,\\
 & \\
\frac{(1-t_i)^{n-j+1} (1-t_{i-j})\prod_{k=1}^{j-1}(t_i-t_{i-k})}{(1-t_{i-1})^{n-j+2}\prod_{k=2}^{j}(t_{i-1}-t_{i-k})}, & j=\frac{n+2}{2},\ldots,n; ~i=j+1,\ldots, n+1,
\end{array} \right. \eqno(3.10)
$$
$$
\widetilde m_{i,j}=\left\{ \begin{array}{ll}
\frac{\frac{n}{2}+i-1}{i-1}t_j, & \quad i=2,\ldots,\frac{n}{2}; ~j=1,\ldots,i-1,\\
\frac{2t_j}{\prod_{k=1}^j (1-t_k)} & \quad i=\frac{n+2}{2}; ~j=1,\ldots,\frac{n}{2},\\
\frac{n-i+2}{\frac{n}{2}+n-i+2}\frac{1}{1-t_j} & \quad i=\frac{n+6}{2},\ldots,n+1; j=1,\ldots,i-\frac{n+4}{2},\\
\frac{n-i+2}{\frac{n}{2}+n-i+2}\frac{t_j}{1-t_j} & \quad i=\frac{n+4}{2},\ldots,n+1; j=i-\frac{n+2}{2},\ldots,i-1,
\end{array} \right. \eqno(3.11)
$$
and
$$
p_{i,i}=\left\{ \begin{array}{ll}
{\frac{n}{2}+i-1 \choose i-1} (1-t_i)^{\frac{n+2}{2}}\prod_{k<i}(t_i-t_k), & \quad i=1,\ldots,\frac{n}{2},\\
 & \\
{\frac{n}{2}+n-i+1\choose n-i+1} \frac{(1-t_i)^{n-i+1} \prod_{k<i}(t_i-t_k)}{\prod_{k=1}^{i-1}(1-t_k)}, & \quad i=\frac{n+2}{2}, \ldots, n+1. \quad \Box
\end{array} \right. \eqno(3.12)
$$

Moreover, by using the same arguments of \cite{MP98}, it can be seen that
this factorization is unique among factorizations of this type, that
is to say, factorizations in which the matrices involved have the
properties shown by formulae (3.4)-(3.6).

\medskip
Let us observe that the formulae obtained in the proof of Theorem
3.2 for the minors of $A$ with $j$ initial consecutive columns and
$j$ consecutive rows, and for the minors of $A^T$ with $j$ initial
consecutive columns and $j$ consecutive rows show that they are not
zero, and so the complete Neville elimination of $A$ can be
performed without row and column exchanges. Looking at equations
(3.7)-(3.12) it is easily seen that $m_{i,j}$, $\widetilde m_{i,j}$ and
$p_{i,i}$ are positive. Therefore, taking into account Theorem 2.1,
this confirms that the matrix $A$ is strictly totally positive.

\medskip
{\bf Theorem 3.3.} Let $A=(a_{i,j})_{1\le i,j \leq n+1}$ be a
SB--Vandermonde matrix  whose nodes satisfy $0 < t_1 < t_2
<\ldots < t_n <t_{n+1} <1$. Then $A$ admits a factorization in
the form
$$A=F_{n}F_{n-1} \cdots F_1 D G_1 \cdots G_{n-1}G_{n}$$
where $F_i$ are lower triangular bidiagonal matrices, $G_i$ are upper triangular ($i=1,\ldots,n$), and $D$ is a diagonal matrix.

{\bf Proof.} The matrix $A$ is a strictly totally matrix \cite{DP06} and therefore, , by Theorem 2.1, the complete Neville elimination of
$A$ can be performed without row and column exchanges providing the
following factorization of $A$ (see \cite{GP96}):
$$A=F_{n}F_{n-1} \cdots F_1 D G_1 \cdots G_{n-1}G_{n},$$
where $F_i$ ($1\leq i \leq n$) are bidiagonal matrices of the form
$$F_i=\left( \begin{array}{cccccccc}
1 & & & & & & & \\
0 & 1 & & & & & & \\
& \ddots & \ddots & & & & & \\
& & 0 & 1 & & & & \\
& & & m_{i+1,1} & 1 & & & \\
& & & & m_{i+2,2} & 1 & & \\
& & & & & \ddots & \ddots & \\
& & & & & & m_{n,n-i} & 1
\end{array}\right),  \eqno(3.9)
$$
$G^T_i$ ($1\le i\le n$) are bidiagonal matrices of the form
$$G_i^T=\left(\begin{array}{cccccccc}
1 & & & & & & & \\
0 & 1 & & & & & & \\
& \ddots & \ddots & & & & & \\
& & 0 & 1 & & & & \\
& & & \widetilde m_{i+1,1} & 1 & & & \\
& & & & \widetilde m_{i+2,2} & 1 & & \\
& & & & & \ddots & \ddots & \\
& & & & & & \widetilde m_{n,n-i} & 1
\end{array}\right),  \eqno(3.10)
$$
and $D$ is the  diagonal matrix
$$
D =\textnormal{diag}\{p_{1,1},p_{2,2},\ldots,p_{n+1,n+1}\}.
$$
The expressions of the multipliers $m_{i,j}$ $(1\leq i,j\leq n+1)$ of the Neville elimination of $A$, the multipliers $\widetilde m_{i,j}$ $(1\leq i,j\leq n+1)$ of the Neville elimination of $A^T$, and the diagonal pivots $p_{i,i}$ $(1\leq i,\leq n+1)$ of the Neville elimination of $A$ are also in this case the ones given by Eq. (3.7) and Eq. (3.10), Eq. (3.8) and Eq. (3.11),  and Eq. (3.9) and Eq. (3.12), respectively. $\Box$

It must be observed that the matrices $F_i$ and $G_i$ ($i=1,\ldots, n$) that appear in the bidiagonal factorization of $A$ are not the same bidiagonal matrices that appear in the bidiagonal factorization of $A^{-1}$ , nor their inverses (see Theorem 3.2 and Theorem 3.3). The multipliers of the Neville elimination of $A$ and $A^T$ give us the bidiagonal factorization of $A$ and $A^{-1}$, but obtaining the bidiagonal factorization of $A$ from the bidiagonal factorization of $A^{-1}$ (or vice versa) is not straightforward. The structure of the bidiagonal matrices that appear in both factorizations is not preserved by the inversion, that is, in general, $F_i^{-1}$ and $G_i^{-1}$ $(1\leq i,j\leq n)$ are not bidiagonal matrices. See \cite{GP96} for a more detailed explanation.

\section{The algorithm} %%%%%%%%%%%%%%%%%%%%%%%%%%%%%%%%%%%%%%%%%%%%%%%%%%%%%%%%%%%%%%%%%%%%%%%%%%%%%%%%%%%%%%%%%%%%%%%%%%%%%%%

In this section we present a fast and accurate algorithm for computing $\mathcal{BD}(A)$ for a totally positive SB--Vandermonde matrix $A$. Let us point out here that given $A$ the  matrix $\mathcal{BD}(A)$ represents {\it both} the bidiagonal decomposition of $A$, and that of its inverse $A^{-1}$ (see Theorem 3.2 and Theorem 3.3).

The algorithm will compute the multipliers $m_{ij}$ of the Neville elimination of $A$, the multipliers $\widetilde m_{ij}$ of the Neville elimination of $A^T$ and the diagonal pivots $p_{ii}$ of the Neville elimination of $A$, which are the entries of the matrix $\mathcal{BD}(A)$.

%In this section we present a fast algorithm for solving the linear system $Ax=b$ corresponding to the following interpolation problem in the Said-Ball basis: given the interpolation nodes $\{t_i: ~i=1, \ldots, n+1\}$ satisfying $0<t_1<t_2< \ldots<t_{n+1}$ and the interpolation data $\{b_i: ~i=1, \ldots, n+1\}$ find the polynomial
%$$
%p(t)=\sum_{k=0}^n a_k s_k^n(t)
%$$
%such that $p(t_i)=b_i$ for $i=1,\ldots,n+1$. In order to solve this linear system whose coefficient matrix $A$ is the $(n+1)\times(n+1)$ SB--Vandermonde matrix introduced in Section 3, we use Theorem 3.2 for obtaining
%$$x=A^{-1}b=G_1G_2\cdots G_{n}D^{-1}F_{n}F_{n-1}\cdots F_1b.$$
%Since $F_i$ and $G_i$ ($i=1,\ldots,n+1$) are bidiagonal matrices and
%$D^{-1}$ is a diagonal matrix, it is clear that the computational
%complexity of computing the whole product from right to left is
%$O(n^2)$. It remains to see that the construction of the matrices
%$F_i$, $G_i$ and $D^{-1}$ can be carried out with a computational
%complexity of $O(n^2)$.

We include here the algorithm for the case in which $n$ is an odd number, the algorithm for the even case being analogous.

The algorithm for  computing the $m_{i,j}$ given by Eq. (3.7) is:

{\tt for} $i=2:n+1$

\quad $m_{i,1}=\frac{(1-t_i)^{\frac{n+1}{2}}}{(1-t_{i-1})^{\frac{n+1}{2}}}$

\quad {\tt for} $j=1:\min(i-2,\frac{n-1}{2})$

\quad \quad $m_{i,j+1}=\frac{t_i-t_{i-j}}{t_{i-1}-t_{i-j-1}}\cdot m_{i,j}$

\quad {\tt end}

{\tt end}

{\tt for} $i=\frac{n+5}{2}:n+1$

\quad $m_{i,\frac{n+3}{2}}=\frac{(1-t_{i-\frac{n+3}{2}})(t_i-t_{i-\frac{n+1}{2}})}{(1-t_i)(t_{i-1}-t_{i-\frac{n+3}{2}})} \cdot m_{i,\frac{n+1}{2}}$

\quad {\tt for} $j=\frac{n+3}{2}:i-2$

\quad \quad $m_{i,j+1}=\frac{(1-t_{i-1})(1-t_{i-j-1})(t_i-t_{i-j})}{(1-t_i)(1-t_{i-j})(t_{i-1}-t_{i-j-1})} \cdot m_{i,j}$

\quad {\tt end}

{\tt end}

The algorithm for the computation of the $\widetilde m_{i,j}$ given by Eq. (3.8) is:

{\tt for} $i=2:\frac{n+1}{2}$

\quad $aux=\frac{\frac{n-1}{2}+i-1}{i-1}$

\quad for $j=1:i-1$

\quad \quad $\widetilde m_{i,j}=aux\cdot t_j $

\quad {\tt end}

{\tt end}

$\widetilde m_{\frac{n+3}{2},1}=\frac{t_1}{1-t_1}$

{\tt for} $j=1:\frac{n-1}{2}$

\quad $\widetilde m_{\frac{n+3}{2},j+1}= \frac{t_{j+1}}{t_j(1-t_{j+1})}\cdot \widetilde m_{\frac{n+3}{2},j}$

{\tt end}

{\tt for} $i=\frac{n+5}{2}:n+1$

\quad $aux=\frac{n-i+2}{\frac{n-1}{2}+n-i+2}$

\quad {\tt for} $j=1:i-\frac{n+3}{2}$

\quad \quad $int=\frac{1}{1-t_j}$

\quad \quad $\widetilde m_{i,j}=aux \cdot int$

\quad {\tt end}

\quad {\tt for} $j=i-\frac{n+1}{2}:i-1$

\quad \quad $int=\frac{t_j}{1-t_j}$

\quad \quad $\widetilde m_{i,j}=aux \cdot int$

\quad {\tt end}

{\tt end}

The algorithm for computing the diagonal pivots $p_{i,i}$  given by Eq. (3.9) is:

$q=1$

$p_{1,1}=(1-t_1)^{\frac{n+1}{2}}$

{\tt for} $i=1:\frac{n-1}{2}$

\quad $q=\frac{\frac{n-1}{2}+i}{i}\cdot q$

\quad $aux=1$

\quad {\tt for} $k=1:i$

\quad \quad $aux=(t_{i+1}-t_k)\cdot aux$

\quad {\tt end}

\quad $p_{i+1,i+1}=q \cdot (1-t_{i+1})^{\frac{n+1}{2}}\cdot aux$

{\tt end}

$aux=1$

{\tt for} $k=1:\frac{n+1}{2}$

\quad $aux=(1-t_k)\cdot aux$

{\tt end}

$q=\frac{q}{aux}$

$aux=1$

{\tt for} $k=1:\frac{n+1}{2}$

\quad $aux=(t_{\frac{n+3}{2}}-t_k) \cdot aux$

{\tt end}

$p_{\frac{n+3}{2},\frac{n+3}{2}}=q \cdot (1-t_{\frac{n+3}{2}})^{n-\frac{n+1}{2}}$

{\tt for} $i=\frac{n+3}{2}:n$

\quad $q=\frac{n-i+1}{\frac{n-1}{2}+n-i+1}\cdot \frac{1}{1-t_i} q$

\quad $aux=1$

\quad {\tt for} $k=1:i$

\quad \quad $aux=(t_{i+1}-t_k) \cdot aux$

\quad {\tt end}

\quad $p_{i+1,i+1}=q \cdot (1-t_{i+1})^{n-i} \cdot aux$

{\tt end}

Looking at this algorithm is enough to conclude that:

\begin{itemize}
\item[-] The computational complexity of the computation of $m_{ij}$, $\widetilde m_{ij}$ and $p{ii}$, i.e. of the computation of  $\mathcal{BD}(A)$ is $O(n^2)$.

\item[-] The algorithm has high relative accuracy because it only involves arithmetic operations that avoid inaccurate cancellation.

\item[-] The algorithm does not construct the SB--Vandermonde matrix, it only works with the nodes $\{t_i\}_{1\leq i\leq n+1}$.
\end{itemize}

As for the even case, the properties of the algorithm are exactly the same.

\section{Accurate computations with SB--Vandermonde matrices}

In this section algorithms for solving linear systems and for eigenvalue computation are presented for the case of a totally positive SB--Vandermonde matrix $A$. The algorithms are both accurate and efficient and are based on the algorithm presented in Section 4 for computing $\mathcal{BD}(A)$.

Let us observe here that, of course, one could try to solve these problems by using standard algorithms. However
 the solution provided by them will generally be less accurate since
SB--Vandermonde matrices are ill conditioned (see the numerical experiments in Section 6) and
these algorithms can suffer from inaccurate cancellation, since they do not take into account the structure of
the matrix, which is crucial in our approach.

\subsection{Linear system solving}

Let $Ax=b$ be a linear system whose coefficient matrix $A$ is a SB--Vandermonde matrix of order $n+1$ generated by the nodes $\{t_i\}_{1\leq i\leq n+1}$, where $0<t_1<\ldots<t_{n+1}<1$.

The following algorithm solves $Ax=b$ in a fast way.

INPUT: The nodes $\{t_i\}_{1\leq i\leq n+1}$ and the data vector $b\in {\bf R}^{n+1}$.

OUTPUT: The solution vector $x\in {\bf R}^{n+1}$.

\begin{itemize}
\item[-] {\it Step 1:} Computation of $\mathcal{BD}(A)$ by using the algorithm introduced in Section 4.

\item[-] {\it Step 2:} Computation of $$x=A^{-1}b=G_1G_2\cdots G_{n}D^{-1}F_{n}F_{n-1}\cdots F_1b.$$
\end{itemize}

Step 2 can be carried out by using the algorithm {\tt TNSolve} of P. Koev \cite{KOEV}. Given the bidiagonal factorization of the matrix $A$, {\tt TNSolve} solves $Ax=b$ by computing the above matrix product.

Although $\mathcal{BD}(A)$ is computed with high relative accuracy, the accuracy of the solution vector will generally depend on the data vector $b$ \cite{MM07}.

Taking into account that, as we have shown in Section 4, the computational cost of Step 1 is of $O(n^2)$ arithmetic operations, and the cost of computing whole product in Step 2 (from right to left) is also of  $O(n^2)$ arithmetic operations, the computational complexity of the algorithm for solving $Ax=b$ is $O(n^2)$.

\subsection{Eigenvalue computation}

Let $A$ be a SB--Vandermonde matrix of order $n+1$ generated by the nodes $\{t_i\}_{1\leq i\leq n+1}$, where $0<t_1<\ldots<t_{n+1}<1$. The following algorithm computes accurately the eigenvalues of $A$.

INPUT: The nodes $\{t_i\}_{1\leq i\leq n+1}$.

OUTPUT: A vector $x\in {\bf R}^{n+1}$ containing the eigenvalues of $A$.

\begin{itemize}
\item[-] {\it Step 1:} Computation of $\mathcal{BD}(A)$ by using the algorithm introduced in Section 4.

\item[-] {\it Step 2:} Given the result of Step 1, computation of the eigenvalues of $A$ by using the algorithm {\tt TNEigenvalues}.
\end{itemize}

{\tt TNEigenvalues} is an algorithm of P. Koev \cite{KOEV05} which computes accurate eigenvalues of a totally positive matrix starting from its bidiagonal factorization. The computational cost of {\tt TNEigenvalues} is of $O(n^3)$ arithmetic operations (see \cite{KOEV05}) and its implementation in \textsc{Matlab} can be taken from \cite{KOEV}. In this way, as the computational cost of Step 1 is of $O(n^2)$ arithmetic operations, the cost of the whole algorithm is of $O(n^3)$ arithmetic operations.

\section{Numerical experiments}

In this section we present two numerical experiments illustrating the accuracy of the two algorithms we have introduced in the previous section.

{\bf Example 6.1.} Let $\mathcal{S}_{15}$ be the Said--Ball basis of the space of polynomials with degree less than or equal to $15$ in $[0,1]$, and let $A$ be the SB--Vandermonde matrix of order $16$ generated by the following nodes:
{\scriptsize$$
\frac{1}{16}< \frac{1}{13} < \frac{2}{11} < \frac{3}{13} < \frac{1}{4} < \frac{7}{18} < \frac{2}{5} < \frac{4}{9} < \frac{7}{15} < \frac{17}{30} < \frac{15}{26} < \frac{9}{13} < \frac{7}{10} < \frac{8}{11} < \frac{5}{6} < \frac{20}{21}.
$$}
The condition number of $A$ is: $\kappa_2(A)=3.2e+08$. Let us consider the data vector
$$
b=(12,-3,0,1,5, -7,0,2, 21,-4,0,9,-11,6,-8,0)^T.
$$
We compute the exact solution $x_e$ of the linear system $Ax=b$ by using the command {\tt linsolve} of {\it Maple 10} and we use it for comparing the accuracy of the results obtained in {\sc Matlab} by means of:

\begin{enumerate}

\item The algorithm presented in Section 5.1. We will call it {\tt MM}.

\item The algorithm {\tt TNBD} of Plamen Koev \cite{KOEV} that computes $\mathcal{BD}(A)$ without taking into account the structure of $A$.

\item The command $A\backslash b$ of {\sc Matlab}.

\end{enumerate}

In $(2)$, the second stage in the solution of the linear system is the computation of the fast product (from right to left) of the bidiagonal matrices and the vector $b$. It is done in {\sc Matlab} by using the same command as in (1): {\tt TNSolve} of Koev \cite{KOEV}.

%The fast product (from right to left) of the bidiagonal matrices and the vector $b$ is also done in {\sc Matlab} by means of the command {\tt TNSolve} of P. Koev \cite{KOEV} and is the second stage in the solution of the linear system in $(1)$ and $(2)$.

We compute the relative error of a solution $x$ of the linear system $Ax=b$ by means of the formula:
$$
err=\frac{\parallel x-x_e \parallel_2}{\parallel x_e \parallel_2}.
$$
The relative errors of the solutions of $Ax=b$ computed by means of the approaches $(1)$, $(2)$ and $(3)$ are reported in  Table 1.

\begin{table}[h]
\begin{center}
\begin{tabular}{|c|c|c|}
\hline MM & TNBD & $A\backslash b$  \\
\hline  5.1e-16 & 2.2e-09 & 3.9e-10 \\
\hline
\end{tabular}
\end{center}\caption{Relative errors in Example 6.1}
\end{table}

{\bf Example 6.2.} Let $A$ be the SB--Vandermonde matrix of order $16$ considered in Example 6.1. In Table 2 we present the eigenvalues $\lambda_i$ of $A$ and the relative errors obtained when computing them by means of:

\begin{enumerate}
\item The algorithm presented in Section 5.2. We will call it {\tt MM}.

\item The algorithm {\tt TNBD} \cite{KOEV} that computes $\mathcal{BD}(A)$ without taking into account the structure of $A$.

\item The command {\tt eig} from {\sc Matlab}.
\end{enumerate}

In $(2)$, the second stage in the computation of the eigenvalues is done in {\sc Matlab} by using the same command as in (1): {\tt TNEigenvalues} of P. Koev \cite{KOEV}.

The relative error of each computed eigenvalue is obtained by using the eigenvalues computed in {\it Maple 10} with 50-digit arithmetic.
\begin{table}[h]
\begin{center}
\begin{tabular}{|c|c|c|c|}
\hline $\lambda_i$  & MM & TNDB & {\tt eig} \\
\hline $1.0e+00$ & $4.4e-16$ & $1.0e-12$ & $1.8e-15$ \\
\hline $9.4e-01$ & $1.3e-15$ & $2.1e-11$ & $1.2e-15$ \\
\hline $7.0e-01$ & $9.6e-16$ & $2.5e-11$ & $6.4e-16$ \\
\hline $5.2e-01$ & $6.3e-16$ & $1.3e-11$ & $2.1e-16$ \\
\hline $3.1e-01$ & $5.4e-16$ & $7.9e-12$ & $2.7e-15$ \\
\hline $1.4e-01$ & $1.3e-15$ & $1.5e-11$ & $1.3e-15$ \\
\hline $6.0e-02$ & $5.7e-16$ & $1.1e-11$ & $1.1e-15$ \\
\hline $3.0e-02$ & $4.6e-16$ & $6.2e-12$ & $4.6e-16$ \\
\hline $8.6e-03$ & $4.1e-16$ & $4.6e-12$ & $1.3e-14$ \\
\hline $2.6e-03$ & $9.9e-16$ & $1.0e-11$ & $3.7e-14$ \\
\hline $6.1e-04$ & $5.4e-16$ & $2.3e-11$ & $7.2e-14$ \\
\hline $6.2e-05$ & $0$  & $1.0e-11$ & $3.1e-13$ \\
\hline $8.3e-06$ & $4.1e-16$ & $1.8e-11$ & $6.4e-13$ \\
\hline $9.1e-07$ & $1.2e-16$ & $4.6e-11$ & $3.2e-12$ \\
\hline $5.5e-08$ & $2.0e-15$ & $1.2e-10$ & $3.1e-10$ \\
\hline $5.0e-09$ & $3.0e-15$ & $2.3e-09$ & $2.0e-09$ \\
\hline
\end{tabular}
\end{center}\caption{Relative errors in Example 6.2}
\end{table}

The results appearing in Table 1 and Table 2 illustrate the good behaviour of our approach. In particular, the very different results obtained for the approaches (1) and (2) show the importance of computing $\mathcal{BD}(A)$ with high relative accuracy, since in both approaches the second stage is exactly the same.

For this specific matrix $A$ the relative error obtained when computing the matrix $\mathcal{BD}(A)$ by using the algorithm we have presented in Section 4 is
$2.8e-15$, while the relative error obtained when computing it by means of the command {\tt TNBD} is
$6.8e-10$. These relative errors have been computed for each solution $B$ by using
$$
err=\frac{\parallel B-B_e \parallel_2}{\parallel B_e \parallel_2},
$$
where $B_e$ is the exact  $\mathcal{BD}(A)$ computed in {\it Maple 10}.

\section*{Acknowledgements}

This research has been partially supported by Spanish Research Grant
MTM2006-03388 from the Spanish Ministerio de Educaci\'on y Ciencia.

%\label{}

% \appendix

% \section{}
% \label{}

% Bibliographic references with the natbib package:
% Parenthetical: \citep{Bai92} produces (Bailyn 1992).
% Textual: \citet{Bai95} produces Bailyn et al. (1995).
% An affix and part of a reference:
%   \citep[e.g.][Ch. 2]{Bar76}
%   produces (e.g. Barnes et al. 1976, Ch. 2).

\end{document}